\documentclass[11pt]{article}
\usepackage{hyperref}
\usepackage[english]{babel}

\usepackage{graphicx}
\usepackage{amssymb}
\usepackage{amsthm, amsmath,amssymb}
\usepackage[all]{xy}
\usepackage{pstricks,pstricks-add,pst-node,pst-text,pst-3d}
\usepackage{graphicx}
\usepackage{psfrag}
\usepackage{multirow}

\usepackage[applemac]{inputenc}

\newcommand{\mathsym}[1]{{}}
\newcommand{\unicode}[1]{{}}

\newtheorem{theorem}{Theorem}[section]

\theoremstyle{definition}
\newtheorem{definition}[theorem]{Definition}

\newtheorem{exemplo}{Example}

\begin{document}

\title{Four beautiful quadrature rules}

\author{M{\' a}rio M. Gra{\c c}a\thanks{Departamento de Matem\'{a}tica,
Instituto Superior T\'ecnico, Universidade  de Lisboa,
 Av. Rovisco Pais,                 
1049--001 Lisboa, Portugal, 
e-mail: \texttt{\url{mario.meireles.graca@tecnico.ulisboa.pt}}  }
}
 
\maketitle

\begin{abstract}

\noindent
A framework is presented to compute approximations of an integral $I(f)=\displaystyle \int_a^b f(x) dx$ from a pair of  companion rules and its associate rule. We show  that an associate rule is a  weighted mean of two companion rules. In particular, the trapezoidal (T) and Simpson (S) rules are  weighted means  of the companion pairs (L,R) and (T,M) respectively,  with L  the  left  rectangle, R the  right rectangle and M the midpoint rules.  As L,R,T and M   reproduce exactly the number $\pi=\displaystyle \int_0^\pi 2\, \sin^2(x) dx$,  we named them the four \lq{beautiful}\rq \  rules. For this example the geometrical interpretation of the rules suggest possible applications of the transcendental number $\pi$ in  architectural design, justifying the attribute \lq{beautiful}\rq\  given to the mentioned rules. 
 As a complement  we consider other appropriate integrand functions $f$, applying  composite rules in order to obtain good approximations of $\pi$, as shown in the worked numerical examples.   
\end{abstract}



 \section{Introduction}\label{introd}

\noindent
From the antiquity to Newton \cite{newton} and  Gauss  \cite{gauss} times, and from Gauss to nowadays a great amount  of ingenious work to approximate an integral $I(f)=\displaystyle \int_a^b f(x) \,dx$ took place giving rise to a plethora of  quadrature rules \cite{newton}, \cite{whittaker}, \cite{davis}, \cite{krylov}, \cite{abramowitz}, \cite{nikolski}, \cite{stroud}, \cite{krommer}. 
We will restrict the discussion to the simplest quadrature rules,  known as  {\em left  rectangle} (L), {\em right rectangle} (R),  {\em midpoint}  (M) and  {\em trapezoidal}  (T).  We call  them the four \lq{beautiful}\rq\ rules since all exactly reproduce  the number $\pi=\displaystyle \int_0^\pi 2\, \sin^2(x) dx$, as shown in a worked example (see Section \ref{exemplos} and Figures \ref{graal}, \ref{abstract} and \ref{figcathedral}), and the respective geometrical interpretation suggests an application of this transcendental number in architectural design.

\medskip
\noindent
After defining  two main types of rules, the Taylor and interpolatory ones, we assign a  {\em sign} to a rule  (Definition \ref{def2}). From a pair of basic rules \--- with opposite errors and the same {\em degree of precision} \--- hereafter called  {\em companion rules},  we construct another rule, which will be called as {\em associate} rule.
In particular, we show that the trapezoidal rule is the associate rule to the pair (L,R), and the famous Simpson's rule is the associate rule to the pair (M,T). 

\medskip
\noindent
Similar (elementar) approach can be applied to other well known sets of quadrature rules, namely to certain pairs of closed and open Newton--Cotes rules \cite{graca}, but we keep our discussion restricted  to a small set of simpler  rules. 

\medskip
\noindent
The consideration of companion rules enables the construction of a nest of intervals where it lies the exact value $I(f)$. Thus, one automatically obtain error bounds for a triple of rules (X,Y,Z) where X and Y are companion rules and $Z$ is the respective associate rule.

\medskip
\noindent
Under mild assumptions, we show that the expression of an associate rule is a {\em weighted mean}, with positive weights, of its companion rules. Therefore,   the value of an associate rule can be obtained by this mean rather than from the explicit expression of  the associate rule. For instance, denoting by $\tilde X$  the value obtained by a rule $X$, one may compute the value of Simpson's rule $\tilde S$ by a certain weighted mean of $\tilde M$ and $\tilde T$, and consequently to be sure that in the interval $[\tilde M, \tilde T]$, or $[\tilde T, \tilde S]$, lies not only $\tilde S$ but also the value of the integral $I(f)$, whenever $(M,T)$ is a pair of companion rules.

\medskip
\noindent
The main features of our framework are illustrated by some examples given in Section \ref{exemplos}.
 
%========================
\section{Taylor and interpolatory rules}\label{taylor}

\noindent
Firtsly, we recall the notion of a quadrature rule. Let $f$ be a  sufficiently smooth functions on the  interval $[a,b]$. For $k\geq 0$, consider ${\cal P}_k$  the linear space of all polynomials of degree $\leq k$.

\begin{definition} (Quadrature rule \cite{engels}, p. 1)\label{def1}

\noindent
A quadrature rule $Q(f)$ is an approximation of the integral $\displaystyle \int_a^b f(x) dx$ obtained using values of $f$ (and/or its derivatives) on a discrete set of points in $[a,b]$.
\end{definition}

%\noindent
%We first define the two types of quadrature rules to be considered in this work, called respectively Taylor and interpolatory rules.
%
%\noindent
%It is assumed that the given real integrand function $f$ belongs to a set of sufficiently smooth functions on the  interval $[a,b]$. For $k\geq 0$,  we denote by ${\cal P}_k$ the linear space of all polynomials of degree $\leq k$.

\medskip
\noindent
For the sake of simplicity we occasionally denote the functional  $Q(f)$  simply by $Q$.

\medskip
\noindent
The two types of quadrature rules to be considered in this work are Taylor and interpolatory rules, which we briefly review.

\subsection*{Taylor's rules}\label{tay1}
Let $k\geq 0$ be an integer, $x_0\in [a,b]$, and assume that $f\in C^{k+1}[a,b]$. The function $f$ can be written as the sum of the Taylor's polynomial, centered at  $x_0$, and a remainder function $r(x)$, of the form
$$
\begin{array}{ll}
f(x)&=q_k(x)+r(x),\\
\mbox{where}& q_k(x)=f(x_0)+ f'(x_0)\, (x-x_0)+\cdots+ \displaystyle \frac{f^{(k)}}{k!} (x-x_0)^k,\\
\mbox{and}&  \\
r(x)&= \displaystyle \frac{ f^{k+1}(\theta(x))}{(k+1)!}\, (x-x_0)^{k+1}, \quad  \theta(x)\in \, int(x_0,x),
\end{array}
$$
where the symbol $int(...)$ denotes the open interval 
$int(x_0,x)= (min \{ x_0,x\}, max\{x_0,x \}) .$
 Integrating both sides in $f(x)=q_k(x)+r(x)$, we get
\begin{equation}\label{eq2}
I(f)=\displaystyle \int_a^b q_k(x) dx+ \displaystyle \int_a^b \displaystyle \frac{ f^{k+1}(\theta(x))}{(k+1)!}\, (x-x_0)^{k+1} dx,  \quad  \theta(x)\in \, int(x_0,x) \  .
\end{equation}
Assuming that the second term of the sum in \eqref{eq2} is not null, we call  Taylor's rule of order $k$ to
\begin{equation}\label{eq3}
T_k(f) =\int_a^b q_k(t)\, dt.
\end{equation}
Thus, the error of $T_k$ is
\begin{equation}\label{eq4}
E_{T_k}(f) = I(f) -T_k(f)=  \int_a^b \displaystyle \frac{ f^{k+1}(\theta(t))}{(k+1)!}\, (t-x_0)^{k+1} dt .
\end{equation}

\subsubsection*{Degree of a rule}\label{deg}
By construction, a Taylor's rule $T_k(f)$ is exact when  $f$ is a polynomial belonging to ${\cal P}_k$,  that is the respective error is null. This is the reason why  we say that $T_k$ is a rule of {\em degree} of precision at least $k$. When a rule of degree at least $k$ is not exact for a monomial of ${\cal P}_{k+1}$ we say that the rule has degree (of precision) $k$, according to the following definition.

\begin{definition}\label{def2} (Degree of a quadrature formula) (\cite{gautschi}, p. 157)

\noindent
We say that a quadrature rule $Q(f)$ has degree (of precision) $m\geq 0$ if the rule is exact for any polynomial  $p_m\in {\cal P}_m$, but it is not exact for a monomial of degree  $m+1$.
\end{definition}

\subsection*{Interpolatory rules}\label{int}

Given $k+1$ distinct nodes in $[a,b]$, say $x_0<x_1<\ldots<x_k$, consider the table of values $\{ x_i,f(x_i)\}_{i=0}^{i=k}$. It is well known that there exists a unique   polynomial $p_k(t)\in {\cal P}_k$ interpolating the  table.  Assuming that $f\in C^{k+1}[a,b]$, we have
\begin{equation}\label{eq6a}
f(x)=p_k(x)+ R(x)\\
\end{equation}
and there exists $ \theta(x)\in int (x,x_0,x_1,\ldots, x_k)$, such that
$$
R(x)= \displaystyle \frac{f^{(k+1)} (\theta(x)) }{(k+1)!}\, (x-x_0)\,(x-x_1)\cdots (x-x_k)\ .
$$
Integrating both sides of \eqref{eq6a}  we obtain 
$$
I(f)= \int_a^b p_k(x)\, dx+ \int_a^b R(x) dx\ .
$$
Assuming that $\displaystyle  \int_a^b R(x) dx\neq 0$, we call
\begin{equation}\label{eq6d}
A_k(f)=\int_a^b p_k(t) dt
\end{equation}
an interpolatory rule of order $k$. The error of $A_k$ is given by
\begin{equation}\label{eq7}
\begin{array}{ll}
E_{A_k}(f)&=I(f)-A_k(f)\\
\\
&= \displaystyle \int_a^b \displaystyle \frac{f^{(k+1)} (\theta(t)) }{(k+1)!}\, (t-x_0)\,(t-x_1)\cdots (t-x_k)\, dt\ .
\end{array}
\end{equation}
Comparing the error formulas \eqref{eq4} and \eqref{eq7} it is clear that it is easier to deal with the former. This is the reason why it seems more natural to start  discussing  quadrature rules by  the Taylor type ones.

\medskip
\noindent
Whenever $f$ is a polynomial of degree $\leq m$, both rules $T_m(f)$ in \eqref{eq3} and $A_m(f)$ in \eqref{eq6d} are exact. In both cases the expressions of their errors  exhibit the $(m+1)$\--th derivative of $f$, and there exist nonzero constants $\alpha_1,\alpha_2$ such that
\begin{equation}\label{b1}
\begin{array}{l}
E_{T_m}(f)=   \alpha_1 \, \displaystyle \frac{f^{(m+1)} (\xi_1)}{(m+1)!}, \qquad \xi_1\in (a,b)     \\
\\
E_{A_m}(f)=   \alpha_2 \, \displaystyle \frac{f^{(m+1)} (\xi_2)}{(m+1)!}, \qquad \xi_2\in (a,b)\ .     \\   \\
\end{array}
\end{equation}
We are interested in guaranteeing that the  derivatives in an error formula do not change  sign  in  $[a,b]$, which lead us to consider the following assumptions.

%\subsubsection*{Assumptions}

\medskip
\noindent
Consider that $f\in C^{m+1}[a,b]$ and $f^{(m+1)}$ do not change sign in $[a,b]$. That  is, 
\begin{equation}\label{asA}
\hspace{-2cm}\mbox{\bf Assumption A}:\quad sign\left( f^{(m+1)}(x) \right) =1 \,\mbox{or}\,\, -1, \,\, \forall x \in [a,b]
\end{equation}
or 
\begin{equation}\label{asB}
\mbox{\bf Assumption B}:\quad sign\left( f^{(m)}(x)- f^{(m) }(y)  \right) =1 \,\mbox{or}\,\, -1,\,\, \forall x,y,\in [a,b].
\end{equation}
When a function $f$ satisfies Assumption A,  the sign of the errors of $T_m(f)$ or $A_m(f)$ depends only on the sign of the constants $\alpha_1$ or $\alpha_2$ in \eqref{b1}. This justifies the following definition of {\em positive} or {\em negative} rule.

\begin{definition}\label{def3}

\noindent
Let $Q_m(f)$ be a quadrature rule of degree $m \,(m\geq 0)$, satisfying Assumption~A in \eqref{asA}, such that the respective error is of the form
$$
E_{Q_m}(f)=  \alpha\,  \displaystyle f^{(m+1)} (\xi),\quad \alpha\neq 0,  \qquad \xi \in (a,b)\ .
$$
The rule is {\em positive} (resp. negative) if $\alpha>0$ (resp. $\alpha<0$).
\end{definition}

\noindent
A pair of rules of the same degree having errors of opposite signs are particularly interesting. Such rules will be called {\em companion} rules.

\begin{definition}\label{def4}(Companion rules)

\noindent
Two rules $X_m(f)$ and $Y_m(f)$, of the same degree $m$, having opposite signs are called {\em companion rules}.
\end{definition}

\noindent
From the error formulae of two companion rules one can deduce not only a new rule, which we call {\em associate rule}, but also its error expression as follows.

\medskip
\noindent
Without loss of generality assume that the rule $X_m$ is positive and the companion  rule $Y_m$ is negative. Concerning Taylor and interpolatory rules, there are positive integers $d_1>0$ and $d_2>0$ such that the respective error in \eqref{b1} can be written as
\begin{equation}\label{eq10}
E_{X_m}(f)=  \displaystyle \frac{1}{d_1} \, (b-a)^{m+2} \displaystyle f^{(m+1)} (\xi_1), \qquad \xi_1\in (a,b)     \\
\end{equation}
\begin{equation}\label{eq11}
E_{Y_m}(f)=  - \displaystyle \frac{1}{d_2} \, (b-a)^{m+2} \displaystyle f^{(m+1)} (\xi_2), \qquad \xi_2\in (a,b)  \ .   
\end{equation}
Thus, adding \eqref{eq10} and \eqref{eq11} and simplifying we obtain
\begin{equation}\label{eq12}
\begin{array}{ll}
I(f)&=   \displaystyle \frac{d_1\, X_m(f)+ d_2\, Y_m(f)}{d_1+d_2} +\\
\\
&\hspace{0.5cm}+\displaystyle \frac{ (b-a)^{m+2} }{d_1+d_2}\displaystyle \left( f^{(m+1)} (\xi_1)  - f^{(m+1)} (\xi_2) \right), \quad \xi_1,\xi_2\in (a,b) \ .
\end{array}  
\end{equation}

\subsection*{Associate rule}
The sum in \eqref{eq12} shows that we can define a new rule, say $M_m(f)$, by
$$
M_m(f)=   \displaystyle \frac{d_1\, X_m(f)+ d_2\, Y_m(f)}{d_1+d_2},
$$
and the respective error is given by
\begin{equation}\label{eq14}
E_{M_m}(f)=\displaystyle \frac{ (b-a)^{m+2} }{d_1+d_2}\displaystyle \left( f^{(m+1)} (\xi_1)  - f^{(m+1)} (\xi_2) \right), \quad \xi_1,\xi_2\in (a,b) \ .
\end{equation}
   Dividing the natural numbers $d_1$ and $d_2$ by its great common divisor,       we obtain
\begin{equation}\label{eq13a}
M_m(f)=   \displaystyle \frac{c_1\, X_m(f)+ c_2\, Y_m(f)}{c_1+c_2},
\end{equation}
where
\begin{equation}\label{eq13A}
c_1=d_1/gcd(d_1,d_2)\quad\mbox{and}\quad   c_2=d_2/gcd(d_1,d_2)\ .
\end{equation}

\medskip
\noindent
Thus, the rule $M_m(f)$ given by \eqref{eq13a}  is a {\em weighted mean} of the companion rules $X_m(f)$ and $Y_m(f)$. The weights $c_1/(c_1+c_2)$ and $c_2/(c_1+c_2)$ are irreductible positive rational numbers.
 We call $M_m(f)$   the {\em associate rule} to the pair $(X_m(f),Y_m(f))$.  Since such a  weighted mean of the two real numbers  lies between these numbers, the following inequalities holds,
$$
\min\{  X_m(f),Y_m(f) \} \leq M_m(f)\leq  \max\{ X_m(f),Y_m(f) \} \ .
$$

%====================
\section{Basic companion  and  associate rules}\label{prel}
 
  \hspace{5cm}
 {\small
\begin{tabular}{p{7.5cm}} 
\lq\lq The material equipment essential for a student's mathematical laboratory is very simple. Each student should have a copy of Barlow's tables of squares, etc., a copy of Grelle's \lq\lq{Calculating Tables}\rq\rq\  and a seven place table of logarithms.\rq\rq
\vspace{0.2cm}
(Whittaker and Robinson \cite{whittaker}, p. vi).
\end{tabular}
}

\medskip
\noindent
We now show that from the  three basic  Taylor's rules (left rectangle, right rectangle and midpoint rule) we can define pairs of companion rules whose  associate rules are 
of interpolatory type, in particular  the  widely used trapezoidal and Simpson rules. Thus, we  automatically obtain error expressions for these associate rules, as in \eqref{eq14}, which might be compared to the ones available  in the literature. Furthermore, we obtain intervals, containing the exact value of $I(f)$, only taking into account to the properties of the arithmetic and weighted mean of two real numbers. 
 
 \medskip
 \noindent
Denoting by $L$ the left rectangle rule, by $R$ the right rectangle rule (both Taylor's rules of order 0) and by $M$ the midpoint rule (Taylor's of order 1), we show that  $T$ (trapezoidal rule) and $S$ (Simpson rule) are associate rules:
$$
\begin{array}{c}
\mbox{(L, R)} \leadsto T\\
\mbox{(M, T)}\leadsto S,
\end{array}
$$
where the symbol $\leadsto$ is used to refer the associate rule to the pair of companion rules.

\medskip
\noindent
For the sake of completeness we show  that the second order  Taylor's  rule $T2$ is companion to the Simpson's rule and  we obtain
its associate $Q$ as the weighted mean
$$
(T2, S)\leadsto Q=(3\, T2+ 2\, S)/5 \ .
$$
The following table summarizes the main features of the companion and associate rules to be discussed in the next paragraphs.

\begin{center}
 \begin{table}[htbp]
   \tabcolsep=0.10cm
   \hspace{4cm}
 \begin{tabular}{| c |   c | c| c|}\hline
Rule & $Degree$& Associate rule & Degree \\
 \hline
$L^{[+]}$&0 & &  \\
 & &$ T= \displaystyle \frac{L+R}{2} $ &  1\\
 $R^{[-]}$& 1 &  &\\
\hline
 $M^{[+]}$& 1 & $ S= \displaystyle \frac{2\, T+M}{3} $ & 3 \\
 \hline
  $T2^{[+]}$& 3 & $ Q= \displaystyle \frac{2\, T2+3\, S}{5} $ & 3 \\
 \hline
  \end{tabular}
 \caption{ The notation $Q^{[+]}$ refer to the sign of the rule $Q$.\label{t1} }
 \end{table}
 \end{center}
 
 %==========================
 \subsection*{The associate rule of $(L, R)$ is the trapezoidal rule $T$}\label{subA1}
\medskip
\noindent
Let $f\in C^1([a,b])$ and consider the node $x_0=a$ or $x_0=b$. The zero order Taylor expansion of $f$ at $x_0$ gives, respectively,
$$
\begin{array}{ll}
f(x)&=f(a)+ f'(\theta_1(x))\, (x-a), \qquad \theta_1(x)\in int(a,x)\\
f(x)&=f(b)+ f'(\theta_2(x))\, (x-b), \qquad \theta_2(x)\in int(a,x)\ . \\
\end{array}
$$
 Integrating both sides in the above equations we obtain,
 \begin{equation}\label{as1}
I(f)=(b-a)\, f(a)+ \displaystyle \frac{(b-a)^2}{2}\, f'(\xi_1),\qquad \xi_1\in (a,b) \ .
 \end{equation} 
  \begin{equation}\label{as2}
I(f)=(b-a)\, f(b) -  \displaystyle \frac{(b-a)^2}{2}\, f'(\xi_2),\qquad \xi_2\in (a,b) \ .
 \end{equation} 
The left and right rectangle rules are implicitly defined in \eqref{as1} and \eqref{as2}. Namely the respective expressions and errors are:

 \begin{equation}\label{as4}
 \left\{
 \begin{array}{l}
 L(f)= (b-a)\, f(a)\\
 E_L(f)=  \displaystyle \frac{(b-a)^2}{2}\, f'(\xi_1),\qquad \xi_1\in (a,b) \ .
 \end{array}
 \right.
 \end{equation}
  \begin{equation}\label{as5}
  \left\{
 \begin{array}{l}
R(f)= (b-a)\, f(b)\\
 E_R(f)= -  \displaystyle \frac{(b-a)^2}{2}\, f'(\xi_2),\qquad \xi_2\in (a,b) \ .
 \end{array}
 \right.
 \end{equation}
The above two rules are both of degree $m=0$ (see Def. \ref{def2}). Under Assumption A, the rule $L(f)$ is positive while $R(f)$ is negative. The constants $d_1,d_2$ in \eqref{eq10}\--\eqref{eq11} are $d_1=d_2=2$ and then $c_1=c_2=1$. Thus, from \eqref{eq13a} the associate rule to the pair $(L,R)$ is the rule
\begin{equation}\label{as6}
T(f)= \displaystyle \frac{L(f)+R(f)}{2}= \displaystyle \frac{b-a}{2} \left( f(a)+ f(b)  \right) \ .
\end{equation}
 We denoted this associate rule by $T(f)$ (instead of $M(f)$ as in \eqref{eq13a}) since it coincides with the well\--known trapezoidal rule. This is a rule of interpolatory type with two nodes $x_0=a$ and $x_1=b$. Using \eqref{eq14}, the respective error formula is
 \begin{equation}\label{as7}
 E_T(f)= \displaystyle \frac{(b-a)^2}{4}\, \left(f'(\xi_1)- f'(\xi_2)\right),\qquad \xi_1, \xi_2\in (a,b) \ .
 \end{equation}

 %========
  \subsection*{The associate rule to $(M, T)$ is the Simpson's rule $S$}\label{subA2}
Let $f\in C^2([a,b])$, and $x_0=(a+b)/2$. The second order Taylor's rule is obtained from
 $$
 f(x)=f(x_0)+ f+(x_0)\, (x-x_0)+ \displaystyle \frac{f^{(2)} (\theta(x)}{2}\, (x-x_0)^2, \,\, \theta(x)\in int(x_0,x) \ .
 $$
 Integrating both sides of the above equality, we get
 \begin{equation}\label{mid1}
 I(f)= (b-a)\, f\left( \frac{a+b}{2}  \right)+ \displaystyle \frac{(b-a)^3}{24}\, f^{(2)} (\xi_1), \quad \xi_1\in (a,b)\ .
 \end{equation}
 From \eqref{mid1} we recover the  midpoint rule $M(f)$,
  \begin{equation}\label{mid2}
M(f)= (b-a)\, f\left( \frac{a+b}{2}  \right),
 \end{equation}
 whose error is given by the expression
 \begin{equation}\label{mid3}
E_M(f)= + \displaystyle \frac{(b-a)^3}{24}\, f^{(2)} (\xi_1), \quad \xi_1\in (a,b)\ .
 \end{equation}
 Recall that the error of the trapezoidal rule (\cite{suli}, Ch. 7) can be written as
 \begin{equation}\label{mid4}
 E_T(f)= - \displaystyle \frac{(b-a)^3}{12}\, f^{(2)} (\xi_2), \quad \xi_1\in (a,b)\ .
 \end{equation}
  Under Assumption A,  from \eqref{mid1} and \eqref{mid3}, we conclude that $M(f)$ and $T(f)$ are respectively positive and negative and have the same degree $m=1$. So, they are companion rules. The corresponding values $d_1$ and $d_2$ in \eqref{eq12} are $d_1=24$, $d_2=12$, and $gcd(d_1,d_2)=12$. Therefore the values $c_1,c_2$ (see \eqref{eq13A}) are $c_1=2$ and $c_2=1$. Thus, the associate rule to the pair $(M,T)$ is Simpson rule
\begin{equation}\label{mid4}
\begin{array}{ll}
S(f)&=\displaystyle \frac{2\, M(f) +T(f)}{3}\\
\\
&= \displaystyle \frac{b-a}{6}\left[ f(a)+ f(b) + 4\, f\left( \displaystyle \frac{a+b}{2}\right)  \right],
\end{array}
\end{equation}
whose error formula is (see \eqref{eq14})
\begin{equation}\label{mid5}
E_S(f)=  \displaystyle \frac{(b-a)^3}{36}\, \left( f^{(2)} (\xi_1)-f^{(2)} (\xi_2)\right), \quad \xi_1,\xi_2\in (a,b)\ .
\end{equation}
Since   $S(f)$ is positive, by \eqref{mid5} and  \eqref{mid4},  Simpson's rule is a weighted mean of (M,T),  the exact value $I(f)$ belongs to the interval defined by the inequalities
$$
\begin{array}{l}
M(f)\leq   S(f) \leq T(f)
\quad \mbox{or}\quad
T(f)\leq S(f) \leq M(f) \ .
\end{array}
$$
Recall that the second expression in \eqref{mid4} coincides with the famous Simpson's rule.
It is worth mentioning that once computed $M(f)$ by \eqref{mid1} and $T(f)$ by \eqref{as6}, it is more efficient to obtain the value $S(f)$ using the first formula   in \eqref{mid4} rather than the second one.  

\medskip
\noindent
In the following paragraph we need to consider the well\--known error formula for Simpson's rule. Let $h=(b-a)/2$. We have 
\begin{equation}\label{mid7}
\begin{array}{ll}
E_S(f)&=- \displaystyle  \frac{b-a}{180}\, h^4\, f^{(4)} (\xi_2), \quad \xi_2\in(a,b)\\
\\
&= - \displaystyle \frac{(b-a)^5}{180\times 2^4}\, f^{(4)} (\xi_2), \quad \xi_2\in(a,b)\ .
\end{array}
\end{equation}
Thus, under Assumption A in \eqref{asA}, this rule in negative, and the respective degree is $m=3$.
 %===============
  \subsection*{The associate  rule to $(T2, S)$ is a corrected midpoint rule}\label{subA3}
 Let $f\in C^4([a,b])$ and $x_0$ be the midpoint of the interval, i.e., $x_0=(a+b)/2$. Since
 $$
 \begin{array}{ll}
 f(x)&=f(x_0)+ f'(x_0)\, (x-x_0)+ f''(x_0)/2\, (x-x_0)^2+ \\
 &+1/6\, f'''(x_0)\, (x-x_0)^3+ f^{(4)} (\theta(x))/4!\, (x-x_0)^4, \, \theta(x)\in int(x_0,x),
 \end{array}
 $$
and taking into account that $\int_a^b (x-x_0)^j\, dx=0$ for any $j$ odd, the integration of both sides of the above equation, gives
$$
I(f)= (b-a)\, f(x_0)+ \displaystyle \frac{(b-a)^2}{24}\, f''(x_0)+ \displaystyle \frac{(b-a)^5}{1920}\, f^{(4)}(\xi_1), \quad \xi_1\in (a,b)\ .
$$
We denote by $T2(f)$ the  Taylor's rule
$$
\begin{array}{ll}
T2(f)&=(b-a)\, f(x_0)+ \displaystyle \frac{(b-a)^2}{24}\, f''\left(\displaystyle \frac{a+b}{2}\right)\\
&= M(f)+ \displaystyle \frac{(b-a)^2}{24}\, f''\left(\displaystyle \frac{a+b}{2}\right)\ .
\end{array}
$$
whose error expression is
$$
E_{T2}(f)= \displaystyle \frac{(b-a)^5}{1920}\, f^{(4)}(\xi_1), \quad \xi_1\in (a,b)\ .
$$
Recall from \eqref{mid7} that the error for the Simpson rule   is
$$
E_S(f)= - \displaystyle \frac{(b-a)^5}{2\, 880}\, f^{(4)}(\xi_2), \quad \xi_2\in (a,b)\ .
$$
 Under Assumption A in \eqref{asA}, the rule $T2(f)$ is positive while $S(f)$ is negative. The constants $d_1,d_2$ in \eqref{eq10}\--\eqref{eq11} are $d_1=1920$, $d_2=2\, 880$, $d_1+ d_2=4800$, and  $gcd(1920,2880)=960$. So   $c_1=2$, $c_2=3$, and by \eqref{eq13a} the associate rule to  $(T2,S)$ is the rule
\begin{equation}\label{t5}
Q(f)= \displaystyle \frac{2\, T2(f)+ 3 \, S(f)}{5},
\end{equation}
whose error expression has the form
\begin{equation}\label{t6}
E_Q(f)= \displaystyle \frac{(b-a)^5}{4\, 800}\, \left( f^{(4)}(\xi_1) -   f^{(4)}(\xi_1) \right) \quad\xi_1, \xi_2\in (a,b)\ .
\end{equation}
Therefore, under Assumption B in \eqref{asB}, the rule $Q(f)$ is positive. It has the same degree $m=3$ as $T2(f)$ and $S(f)$. However, as the rule $Q$ is the weighted mean \eqref{t5} of $T2$ and $S$, the exact value $I(f)$ lies in
 the interval defined by the inequalities
\begin{equation}\label{t7}
\begin{array}{l}
T2(f)\leq   Q(f) \leq S(f)\quad
\mbox{or}\quad
S(f)\leq Q(f) \leq T2(f) \ .
\end{array}
\end{equation}

\medskip
\noindent
\subsubsection*{Composite rules}

To obtain a composite version of a given simple rule $Q(f)$, one divides the interval $[a,b]$ in $n\geq 2$ parts with uniform length $h=(b-a)/n$, and we apply  $Q$ to each subinterval and finally add the $n$ values obtained.  We denote by $Q_n$ the composite rule corresponding to the application of a simple rule $Q$ to the $n$ subintervals $[a_i,b_i]$, such that $a_i=a+i\, h$, $b_i=a_i+h$, with $i=0,\ldots, (n-1)$.

%=========================
%\newpage
\section{Examples}\label{exemplos}

\begin{figure}[hbt] 
\begin{center}
  \includegraphics[scale=0.420]{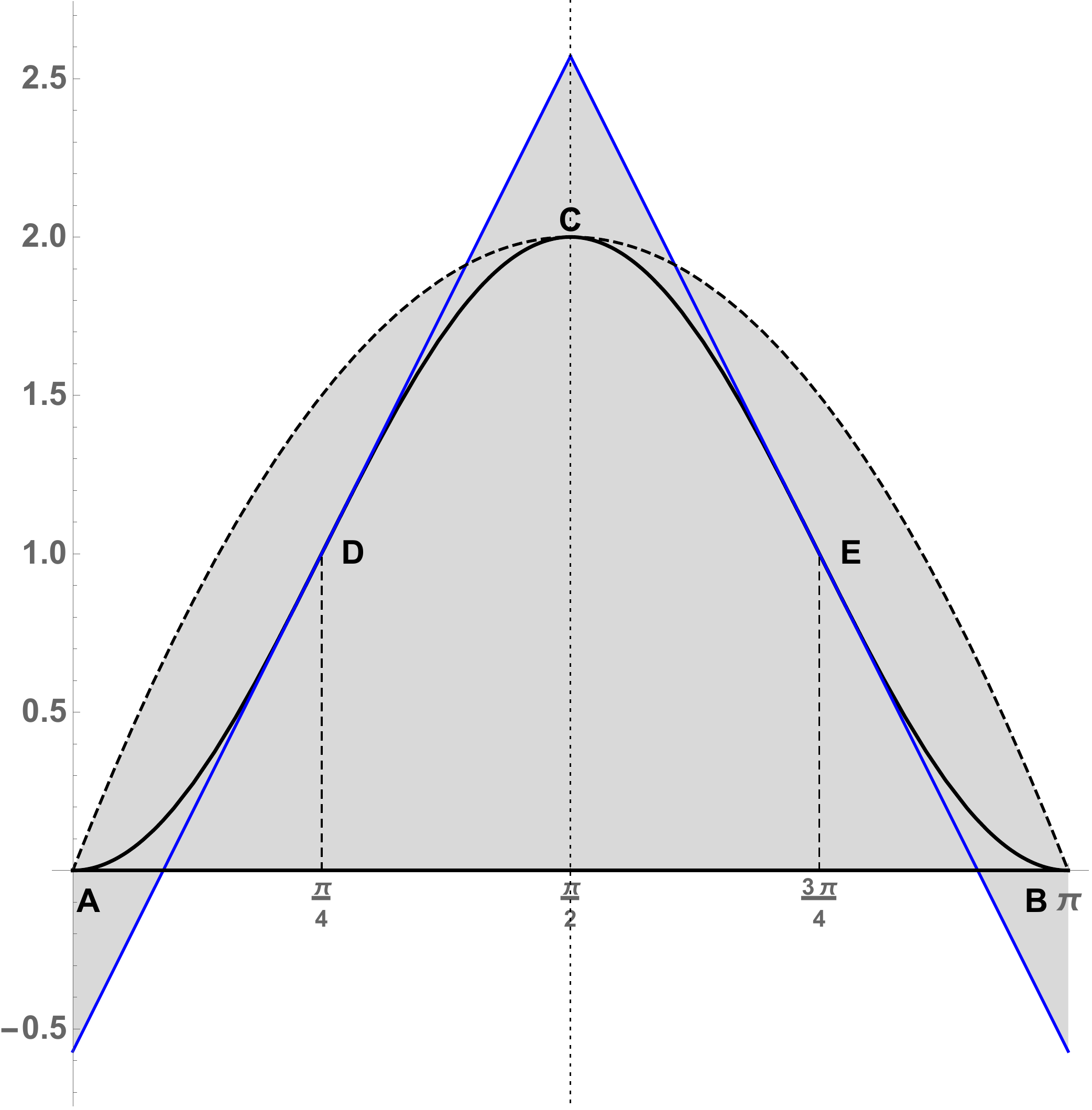} 
 \caption{Four beautiful rules. \label{graal}}
\end{center}
\end{figure}

\begin{figure}[h]
\begin{center} 
 \includegraphics[totalheight=10.cm]{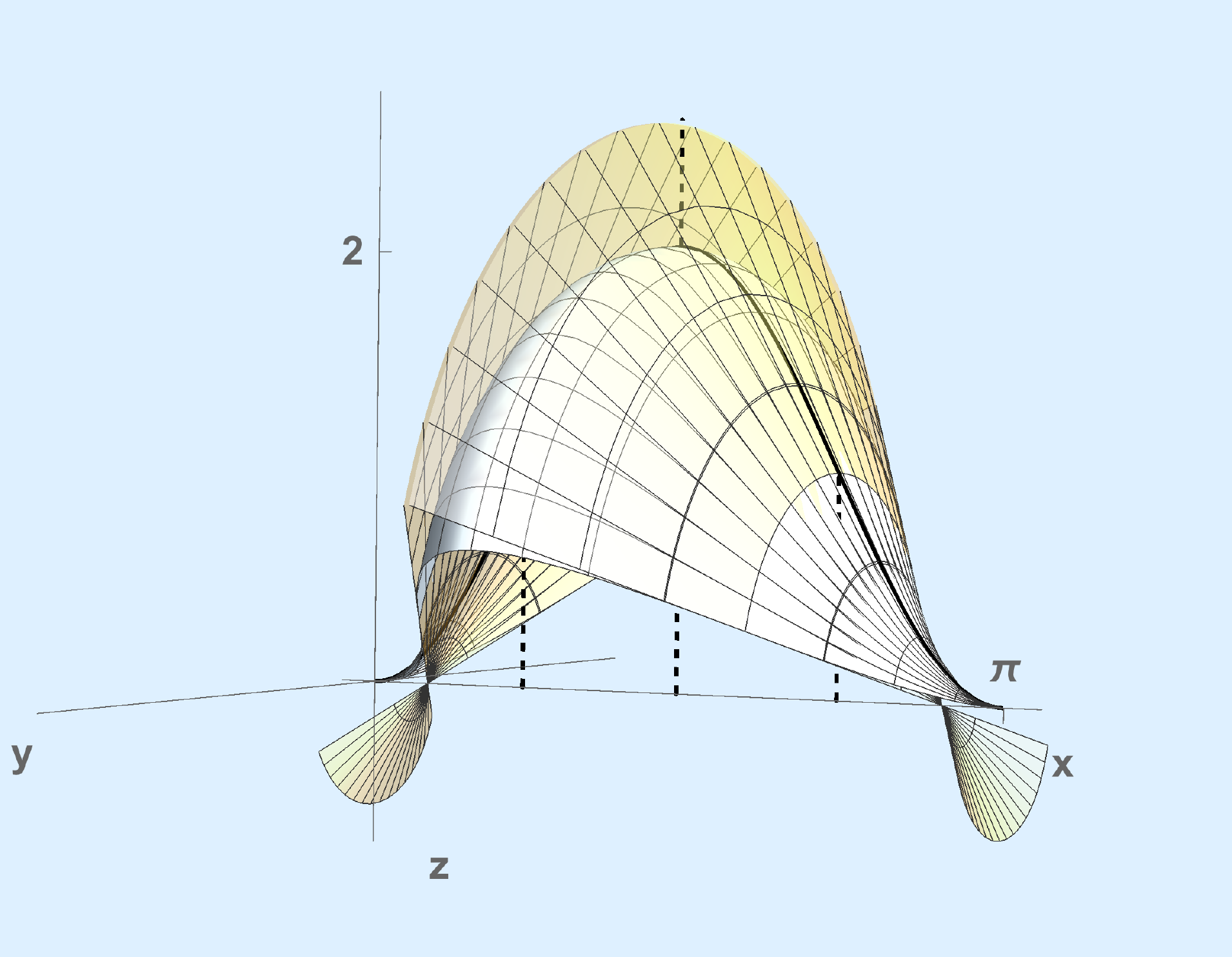}
\caption{\label{figcathedral} Partial revolution around the $x$ axis.  }
\end{center}
\end{figure}

\noindent
The main features of our companion rules and  associate ones can be illustrated by taking as models the following integrals:
$$
\displaystyle \int_0^\pi 2\, \sin^2(x) dx=\displaystyle \int_0^{1/2} \displaystyle \frac{6}{\sqrt{1-x^2}}\, dx=\displaystyle \int_{-1}^{1} \displaystyle \frac{2}{1+x^2}\, dx=\pi\ .
$$
The first integral is particularly interesting since the four composite rules $L_2,R_2,M_2$ and $T2_2$, with $n=2$ subintervals, are exact as shown in  Example \ref{exemplo1}. Thus, the number $\pi$ can be materialized in 2D and 3D pictures as  in Figures \ref{graal} and \ref{figcathedral}, suggesting potential architectural applications.

\medskip
\noindent
The numerical experiments in Examples 2 and 3  were carried out using {\sl Mathematica} \cite{wolfram1}. Commands like
$\verb+Range+$, $\verb+Partition+$ e $\verb+Total+$ are useful to define composite rules (see code given in Example~\ref{exemplo2}).

 \begin{exemplo}\label{exemplo1}
 \end{exemplo}

 \noindent
Let
\begin{equation}\label{ex1}
I(f)=\displaystyle \int_0^\pi 2\, \sin^2(x) dx=\pi\ .
\end{equation}
In Figure \ref{graal} the graph of $f(x)=2\, \sin^2(x)$, passing through the points A,D,C,E, and B is displayed in bold.  The simple rules $L,R,M,T,S$ and $T2$ are useless as approximations of the transcendental number $\pi$, since
$$
L(f)=R(f)=0\quad\mbox{and}\quad 
M(f)=(\pi-0) f(\pi/2)= 2\, \pi \ .
$$
As $T=(L+R)/2$ and $S=(2\, T+M)/3$, we have
$$
T(f)=0\quad \mbox{and}\quad 
S(f)= \displaystyle \frac{2\, \pi}{3}\ .
 $$
Taking into account that $f''(x)=4\, \cos(2\,x)$, the Taylor's rule $T2$ is
$$
T2(f)= M+ \displaystyle \frac{\pi^3}{24}\, f''(\pi/2)= M-\pi^3/6= \pi\, (2-\pi^3/6) \ .
$$
Recall that  $Q_n$ denotes a composite rule with $n$ subintervals of $[a,b]$. The symmetry of the function with respect to the axis $x=\pi/2$ leads to the remarkable property that all the composite rules with $n=2$ subintervals,  $L_2,R_2,M_2$ and the trapezoidal rule $T_2$, produce the exact value of \eqref{ex1}, that is, $\pi$. The subintervals to be considered  are $[0,\pi/2]$ and $[\pi/2, \pi]$. Thus,
$$
\begin{array}{ll}
L_2(f)&=0+\displaystyle \frac{\pi}{2} \times 2=\pi\quad \mbox{(left),}\\
\\
R_2(f)&=\displaystyle \frac{\pi}{2} \times 2+0=\pi\quad \mbox{(right),}\\
\\
T_2(f)&=\displaystyle \frac{\pi+\pi}{2} =\pi \quad \mbox{(trapezoidal),}\\
\\
M_2(f)&=\displaystyle \frac{\pi}{2} \times f\left(\displaystyle \frac{\pi}{4}\right)+\displaystyle \frac{\pi}{2} \times f\left(\displaystyle \frac{\pi}{4}\right)=\pi\quad \mbox{(midpoint)}\ . \\
\end{array}
$$
The simple Simpson's rule ($n=2$ subintervals) is not exact, as suggested by Figure \ref{graal} where the graphic of the  interpolating polynomial passing trough the points A, C, B is displayed using a dashed line. It is obvious that the area delimited by the graph of the interpolating polynomial and the $x$\--axis is not the same as the area under the graph of $f$. However, since the composite rule $M_2$ and $T_2$ are exact, the composite Simpson's rule $S_4$, with $n=4$ subintervals, should be exact because  this rule is a weighted mean of  $M_2$ and $T_2$ given by \eqref{mid4}. In fact,
$$
\begin{array}{l}
S_2(f)= \displaystyle \frac{2\, T(f)+ M(f)}{3}= \displaystyle \frac{ M(f)}{3}= \displaystyle \frac{2\, \pi}{3}\neq \pi,
\end{array}
$$
while
$$
\begin{array}{l}
S_4(f)= \displaystyle \frac{2\, T_2(f)+ M_2(f)}{3}=\displaystyle \frac{3\, \pi}{3}=  \pi \ .
\end{array}
$$
As shown previously that the simple Taylor's rule $T2$ is not exact. However, the composite rule with $n=2$ subintervals is. Let us consider first the subinterval $[0,\pi/2]$. $T2$ represents the area delimited by the  $x$\--axis and  the line segment (passing through the point D  in Fig.\ref{graal}), whose cartesian equation is
$$
y=f(\pi/4)+ f'(\pi/4)\, (x-\pi/4) = 1+ 2\, (x-\pi/4), \quad 0\leq x\leq \pi/2,
$$
 Noting that $f''(\pi/4)=0$, we have
$$
T2(f)= \displaystyle \pi/2\, f(\pi/4)= \displaystyle \frac{\pi}{2}\ .
$$

\begin{figure}[h]
\begin{center} 
 \includegraphics[totalheight=10.cm]{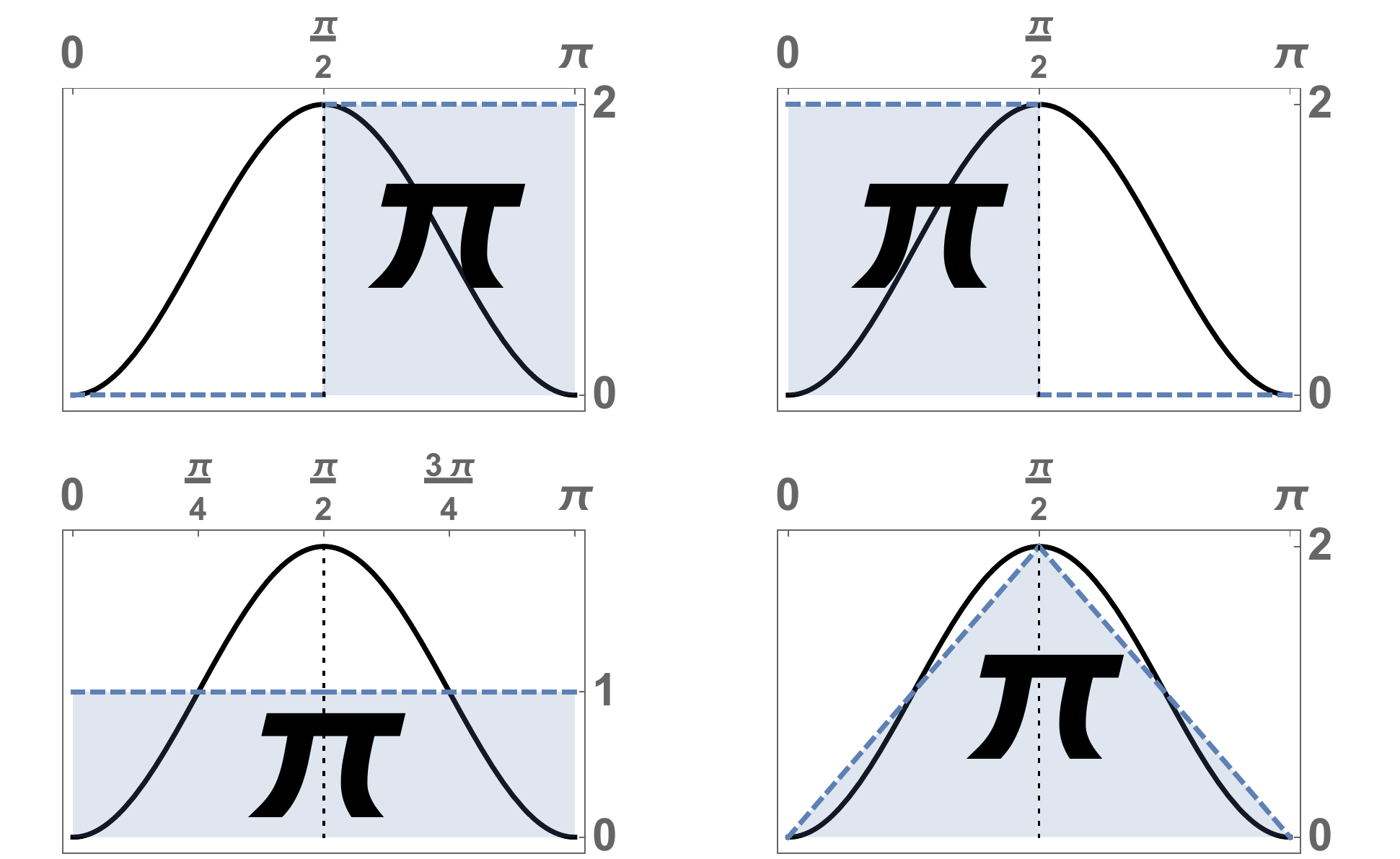}
\caption{\label{abstract} The four beautiful rules $L_2, R_2, M_2$ and $ T_2$. }
\end{center}
\end{figure}

 \noindent
Analogously, in the subinterval $[ \pi/2,\pi   ]$, the value of $T2$ represents the area delimited by the  $x$\--axis and  the line segment (passing through the point E  in Fig.\ref{graal}), whose cartesian equation is
$$
y=f(3\, \pi/4)+ f'(3\, \pi/4)\, (x-3\, \pi/4) = 1- 2\, (x-3\, \pi/4), \quad \pi/2\leq x\leq \pi \ .
$$
As $f''(3\, \pi/4)=0$,  we obtain
$$
T2(f)= \displaystyle \pi/2\, f(3\, \pi/4)= \displaystyle \frac{\pi}{2}\ .
$$
Thus, the composite rule, with $n=2$, subintervals gives
$$
T2_2(f)=\displaystyle \frac{\pi}{2} +\displaystyle \frac{\pi}{2}= \pi \ .
$$
 
\noindent
So, all the above six composite rules, $L_n,R_n,M_n, T_n, S_n,T2_n$, with $n\geq 4$, are exact for the integral \eqref{ex1}. In Figure \ref{abstract} four plane regions whose area is $\pi$ are shown, giving a geometrical interpretation of the rules $L_2$, $R_2$, $M_2$ and $T_2$ respectively.

\medskip
\noindent
 In Figure \ref{figcathedral}  we show  (partial) surfaces generated by the revolution of the graphic  $f$ around the $x$\--axis together with the line segments (cf.  Figure \ref{graal}). This suggests interesting architectural applications of our companion and associate rules.
%=====================
 \begin{exemplo}\label{exemplo2}
 \end{exemplo}
 
 \noindent
 Let 
 $$
 I(f)=\displaystyle \int_0^{1/2} \displaystyle \frac{6}{\sqrt{1-x^2}}\, dx=\pi \ .
 $$
 As the function $f(x)= 6/\sqrt{1-x^2}$ is positive and  its derivative is nonnegative derivative in $[0,1/2]$,  the function $f$ is an increasing function. Therefore, for each $n\geq 1$, the  rules $L_n$ and $R_n$ have errors of opposite signs and
 \begin{equation}\label{ex1a}
 L_n(f)\leq I(f)\leq R_n(f),\quad \forall n\geq 1 \ .
 \end{equation}
 As
 $$
 R_n(f)-L_n(f)=\displaystyle \frac{(b-a)}{n}\, \left( f(b)-f(a) \right)=\displaystyle \frac{1}{2\, n}\, (4 \sqrt{3}-6),
 $$
 we have,
 $$
 \lim_{n\rightarrow \infty} |  R_n(f)-L_n(f)| =0\ .
 $$
 Thus, from \eqref{ex1a} it follows
 $$
  \lim_{n\rightarrow \infty}  R_n(f)=   \lim_{n\rightarrow \infty}  L_n(f)= \lim_{n\rightarrow \infty} =I(f)=\pi\ .
 $$
As the  composite rule $T_n$ is the weighted mean \eqref{as6} of $L_n$ and $R_n$, one  concludes that  $\lim_{n\rightarrow \infty}  T_n(f)=\pi$. Moreover,  as $f$ is an increasing function, we also have for the midpoint rule
$$ 
 L_n(f)\leq M_n(f)\leq R_n(f),\quad \forall n\geq 1\ .
$$
 Consequently, all the rules $L_n,R_n,M_n, T_n$ have $\pi$ as common limit. The same is true  for the Simpson rule $S_n$ since this is the weighted mean \eqref{mid4} between the composite companion rules $M_n$ (midpoint) and  $T_n$ (trapezoidal).

   \begin{figure}[hbt] 
\begin{center}
  \includegraphics[scale=0.65]{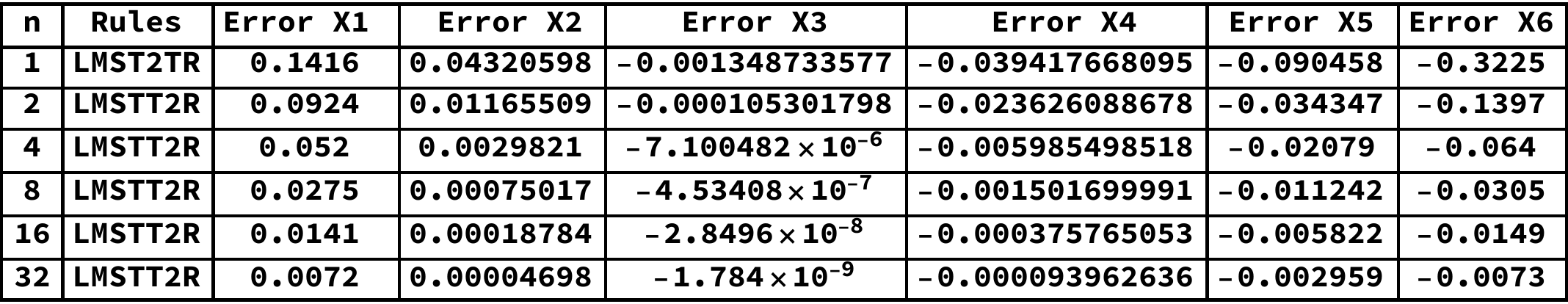} 
 \caption{Composite rules for  $\displaystyle \int_0^{1/2} \displaystyle \frac{6}{\sqrt{1-x^2}}\, dx$ .  \label{figtab}}
\end{center}
\end{figure}

 \medskip
 \noindent
 In Figure \ref{figtab} we compare the errors for each of the previous rules, including  the Taylor' s rule $T2$. For $n=1,2, 2^2, 2^3,\ldots$ we firstly compute the 6-uple $(L_n,R_n,M_n,T_n,S_n,T2_n)$ using the following {\sl Mathematica code}
 \begin{verbatim}
L[f_, a_, b_] := (b - a) f[a];
R[f_, a_, b_] := (b - a) f[b];
M[f_, a_, b_] := (b - a) f[(a + b)/2];
T[f_, a_, b_] := (b - a)/2 (f[a] + f[b]);
S[f_, a_, b_] := (2 M[f, a, b] + T[f, a, b])/3;
T2[f_, a_, b_] := M[f, a, b] + (b - a)^2/24 f''[(a + b)/2];
composite[rule_String, a_, b_, n_] := Block[{h},
   h = (b - a)/n // N;
   Total[Map[ToExpression[rule][f, #[[1]], #[[2]] ] &,
      Partition[Range[a, b, h], 2, 1]]    ] ];
 \end{verbatim}
 
 \medskip
 \noindent 
 The numerical values obtained from the 6 rules are sorted in increasing order. For each $n$, the second column of Figure \ref{figtab}  displays a string reflecting the names of the rule with respect to the order of the computed values. For instance the string LMSTT2R means that the smaller computed value is given by the (composite) left rectangle rule $L_n$ and the greatest one corresponds to the right rectangle rule $R_n$.  The values in the columns  Error Xi refer to the error of the rule in the $i$\--th position of the respective string of rules names.

 \noindent
 All the error columns in Fig. \ref{figtab} show that the respective rule produce monotone sequences. Also, it is clear that $(L_n,R_n)$ and $(M_n,T_n)$ are pairs of companion rules. As the absolute error of $S_n$ is less than the one of $T2_n$,  Simpson's rule performs better than Taylor's $T2$ and all the other rules.
 
 %=====================
% \newpage
 \begin{exemplo}\label{exemplo3}
 \end{exemplo}
 
 \noindent
 Let 
 $$
 I(f)=\displaystyle \int_{-1}^{1} \displaystyle \frac{2}{1+x^2}\, dx=\pi \ .
 $$

 \begin{figure}[h]
\begin{center} 
 \includegraphics[totalheight=3.0cm]{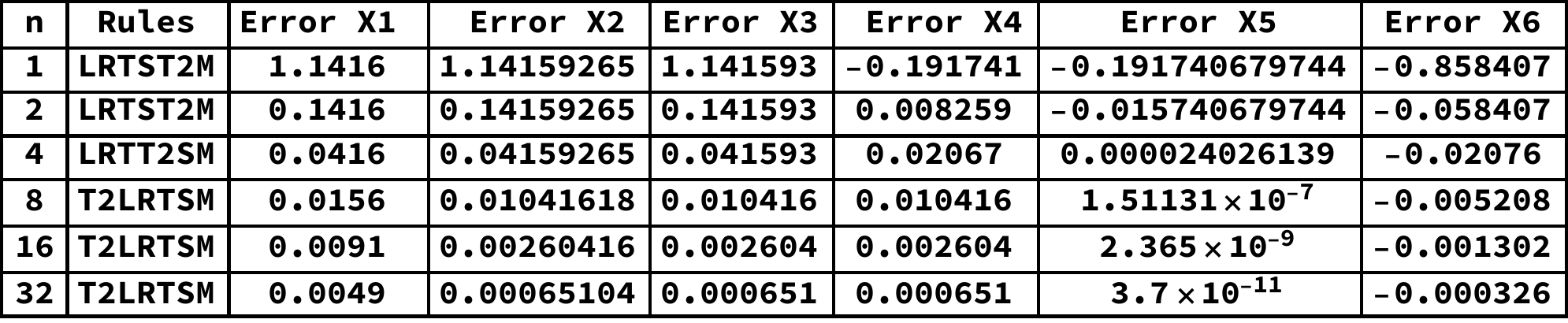}
\caption{\label{tab2} Composite rules for  $\displaystyle \int_{-1}^{1} \displaystyle \frac{2}{1+x^2}\, dx=\pi$ . }
\end{center}
\end{figure}

\medskip
\noindent
In this case $f(x)=2/(1+x^2)$ is a positive function, symmetric with respect to the $y$\--axis. The function increases from $x=-1$ to $x=0$ and decreases from $x=0$ to $x=1$. Both the function and its derivatives are bounded in  $[-1,1]$. As shown in Figure \ref{tab2}, for $1\leq n\leq 4$ the composite rules $L_n$ and $M_n$ (midpoint) have errors of opposite sign, whilst for $8\leq n\leq 32$ the rules $T2_n$ and $M_n$ enjoy the same property.
For all 6 rules their absolute errors seem to approach zero, as expected.

\noindent
One can verify that for $n=2^{10}=1024$ the Simpson's rule gives $\pi\simeq 3.1415926535897932384$, where all digits are correct.
%================

%\newpage

\end{document}